# A Tabu Search Method for the Optimisation of Fluid Power Circuits

**A.M.Connor** and **D.G.Tilley**
Department of Mechanical Engineering, University of Bath.

**Abstract** : This paper describes the development of an efficient algorithm for the optimisation of fluid power circuits. The algorithm is based around the concepts of Tabu search, where different time scale memory cycles are used as a metaheuristic to guide a hill climbing search method out of local optima and locate the globally optimum solution. Results are presented which illustrate the effectiveness of the method on mathematical test functions. In addition to these test functions, some results are presented for real problems in hydraulic circuit design by linking the method to the BATH*fp* dynamic simulation software. In one such example the solutions obtained are compared to those found using simple steady state calculations.

**NOTATION**

| | |
|---|---|
| $e\{k\}$ | Error in term $k$ |
| $T_m$ | Torque produced by motor |
| $D_m$ | Motor displacement |
| $\omega_m$ | Motor speed |
| $\eta_{mm}$ | Motor mechanical efficiency |
| $\eta_{vm}$ | Motor volumetric efficiency |
| $\Delta P_m$ | Pressure drop across motor |
| $Q_m$ | Motor flow rate |
| $Q_p$ | Pump flow rate |
| $Q_{rv}$ | Relief valve flow rate |
| $\omega_p$ | Pump speed |
| $\eta_{vp}$ | Pump volumetric efficiency |

**1 INTRODUCTION**

The manual design of hydraulic circuits is a lengthy process which often involves a large degree of trial and error in the selection and sizing of components which often leads to a compromise in performance in order to achieve results.

The use of numerical optimisation techniques linked with dynamic simulation can lead to better solutions but often require a large amount of time to find the best solution. A significant portion of this time is due to the complex numerical integration required to carry out dynamic simulation of hydraulic circuits.

Whilst some techniques have been developed which require less time to simulate a given circuit, it is more desirable to reduce the number of simulations required by a numerical optimisation method in order to reduce the run time. This is based on previous work (**1,2**) which utilised a Parallel Genetic Algorithm. The large number of simulations required led to unacceptable run time for even very simple circuits.

This paper describes the underlying methodology of the Tabu search that has been tested on a range of mathematical test functions. Results are presented for two numerical test functions and three examples of hydraulic circuit optimisation.

## 2   TABU SEARCH METHOD

Tabu search (**3,4**) is a metaheuristic procedure which is intended to guide optimisation methods to avoid local optima during complex or multi-modal numerical optimisation problems. Tabu search has been successfully applied to a variety of classical and practical problems which include the quadratic assignment problem (**5**), electronic circuit design (**6**) and the balancing of hydraulic turbines (**7**).

The Tabu search method utilises a number of flexible memory cycles, each with different associated time scales, to allow search information to be exploited more thoroughly than by rigid memory or memoryless systems. A brief introduction to the concept of Tabu search will be given here, but more information is contained in (**3,4,8**).

### 2.1 Short Term Memory

The most simple implementation of a Tabu search is based around the use of a hill climbing algorithm. Once the method has located a locally optimal solution, the short term memory is used to force the search out from this optima in a different direction. The short term memory constitutes a form of aggressive search that seeks to always make the best allowable move from any given position. Short term memory is implemented in the form of tabu restrictions which prevent the search from cycling around a given optimal or sub-optimal position. Essentially, this restriction applies to a limited set of previously visited positions which is continuously updated on a "first in, first out" basis as the search progresses.

### 2.2 Aspiration Levels

The use of aspiration levels allows for the tabu status of a given move or position to be temporarily overridden if the move is sufficiently desirable. This is normally implemented in the form that a tabu restriction can be overridden only if the current position is not tabu. This implementation prevents cycling around any given position, but allows the search to explore new areas by revisiting a previously visited position and leave it in a new direction.

### 2.3 Additional Memory Cycles

Intermediate term memory is normally used to intensify the search whilst long term memory is used to diversify the search by refreshing the current search point. The actual implementation of each of these memory cycles will be discussed in sections 2.4.3 and 2.4.4.

**2.4 Specific Implementation**

The simple description of the Tabu search method outlined in the previous sections contains sufficient information for the specific implementation to now be described.

**2.4.1 Hill Climbing Algorithm**

The underlying hill climbing algorithm used in this work is based upon the method developed by Hooke and Jeeves (**9**). This method consists of two stages. The first stage carries out an initial exploration around a given base point. Once this exploration has been carried out and a new point found, the search is extended along the same vector by a factor $k$. This is known as a pattern move. If the pattern move locates a solution with a better objective value than the exploration point, then this point is used as a new base point and the search is repeated. Otherwise, the search is repeated using the exploration point as a new base point. Figure 1 shows a typical exploration and pattern move.

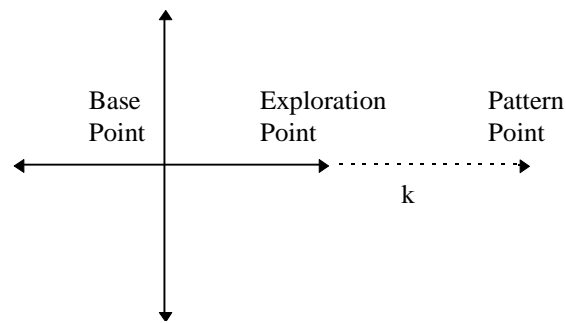

Fig. 1 Hooke–Jeeves search

In the Hooke and Jeeves search, this process is repeated until no further improvement is found. When this occurs, the step size is reduced and the search continued. This is repeated until the step size falls below a given value. In the Tabu search, this search strategy is somewhat modified. Firstly, the selection of the new base point is adjusted so that the search always selects the best available move. When no improvement can be found, the best available move is defined as that where the increase in objective function value is the smallest. Because of this, a new control algorithm has been developed which reduces the step size as the search progresses and eventually terminates the search. This will be described in section 2.4.5

**2.4.2 Tabu Restrictions**

The short term memory cycle has been implemented as a list of Tabu restrictions. The Tabu list contains the parameter values of the last 'n' accepted solutions. During the search, before the objective function of a proposed point is evaluated the parameter set is compared to those contained in the Tabu list. If the position has been recently visited then the objective function is not evaluated and the search classes the move as Tabu.

### 2.4.3 Search Intensification

Any intensification of the search is based around the concept of the intermediate memory cycle. In this case, the intermediate memory is very similar to the short term memory in that it is a list of previously accepted solutions. However, the solutions contained in the intermediate memory list are the previous 'm' best solutions. As a new best solution is found, this is placed into the list and the oldest solution removed. This is termed intermediate memory as the time scale of replacement is much longer than for short term memory.

### 2.4.4 Search Diversification

In many Tabu search applications, search diversification is achieved by the use of long term memory cycles or just by randomly refreshing the current base point when certain conditions apply.

In this case the search is randomly refreshed when the number of trial moves made without an improved solution being found reaches a certain value specified by the user.

### 2.4.5 Search Control

The control algorithm underlying the Tabu search is relatively complex, as it has to ensure that both intensification and diversification occur at appropriate times. In addition to this, the step size must be periodically reduced to ensure that the search is carried out to a reasonably fine detail. The final necessity is that the search must be terminated.

The following psuedocode shows how such control is achieved.

```
start search
control=0
while(terminate!=TRUE)
{
   while(control<END)
   {
    if(control=intense)
        intensify
      if(control=diverse)
        diversify
      search
      if(new best solution)
        control=0
      else
        control=control+1
   }
   reduce step size
   control=0
   if(step size<min step size)
      terminate=TRUE
}
```

At the start of the search, the control flag is set to zero. Whenever a solution is discovered that is better than any previously found solution the control flag is reset to zero. As the search progresses, the control flag is increased by one for each search iteration. When the control flag is equal to certain values then intensification and diversification of the search are carried out. If this still results in no improvement, then when the control flag equals the end of cycle value then it is reset to zero and the step size reduced. The search then repeats until the step sizes are reduced to the stated minimum values when it is terminated.

## 3 TEST FUNCTIONS

The Tabu search has been tested on a variety of mathematical test functions. Two of these functions are outlined here. For each case, the method was run ten times from a random start point and the final solution noted. The number of objective function evaluations for each run were also noted.

### 3.1 Rastrigin Function

The Rastrigin function is a complex two variable problem with one unique minimum and approximately 50 local minima in the permissible search region. The objective function value is calculated using equation 1.

$$obfn = x^2 + y^2 - \cos(18x) - \cos(18y) \tag{1}$$

The search variables are limited to the region enclosed by the bounds $-1.0 \leq x,y \leq 1.0$. The minimum of the function is located at $(x,y) = 0,0$. The Tabu search located the minimum value of the function for all runs on this problem. The average number of objective function evaluations was 1184, given that the initial step size was 0.5 and the minimum step size was 0.0001.

### 3.2 Schwefel Function

The Schwefel function is a very complex ten parameter function with a unique minimum and very large numbers of local minima in the search region. Previous work (**9**) has suggested that methods capable of optimising the Schwefel function are likely to be suitable for optimising hydraulic circuits. The objective function is calculated using equation 2.

$$obfn = \sum_{i=1}^{10} x_i \sin\left(\sqrt{abs(x_i)}\right) \tag{2}$$

The search variables are limited to the region enclosed by the bounds $-500 \leq x_i \leq 500$. The minimum of the function is located at $x_i = 420.9687$. The Tabu search located the minimum value of the function for seven of the ten runs on this problem. The average number of objective function evaluations was 15,159, given that the initial step size was 100 and the minimum step size was 0.01. For the three cases where the optimum solution was not found, the objective function value was within 85% of the function minimum.

### 3.4 Implications of Results

The results obtained from the mathematical test function can be used to assess the potential of Tabu search applied to real problems in hydraulic circuit design. On the

Rastrigin function, the Tabu search consistently located the global optimum with a relatively small number of objective function evaluations. On the Schwefel function the method did not always locate the global optimum. However, the aggressive nature of the Tabu search implies that less objective function evaluations are required than when using alternative optimisation strategies (1,2). As hydraulic circuit simulation is often complex and time consuming, then this aggressive search may lead to the discovery of acceptable solutions to a problem with in a relatively small number of simulations.

**4 HYDRAULIC CIRCUIT OPTIMISATION**

The Tabu search has been used to optimise three hydraulic circuits. The first of these, shown in Figure 2, is a simple hydrostatic transmission. The second circuit, shown in Figure 3, involves a single pump driving two motors at different speeds. The final circuit, shown in Figure 4, consists of a linear actuator that is controlled using proportional feedback in order to generate a desired extend and retract motion profile.

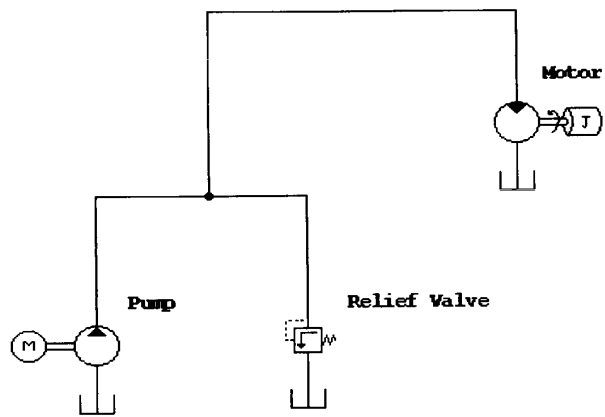

Fig. 2 Hydrostatic transmission

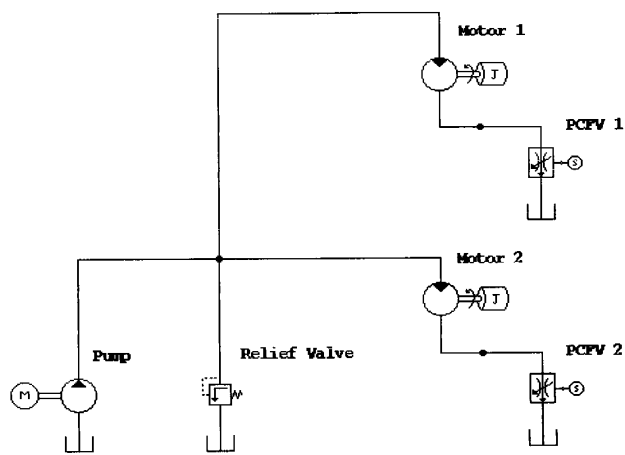

Fig. 3 Two-motor circuit

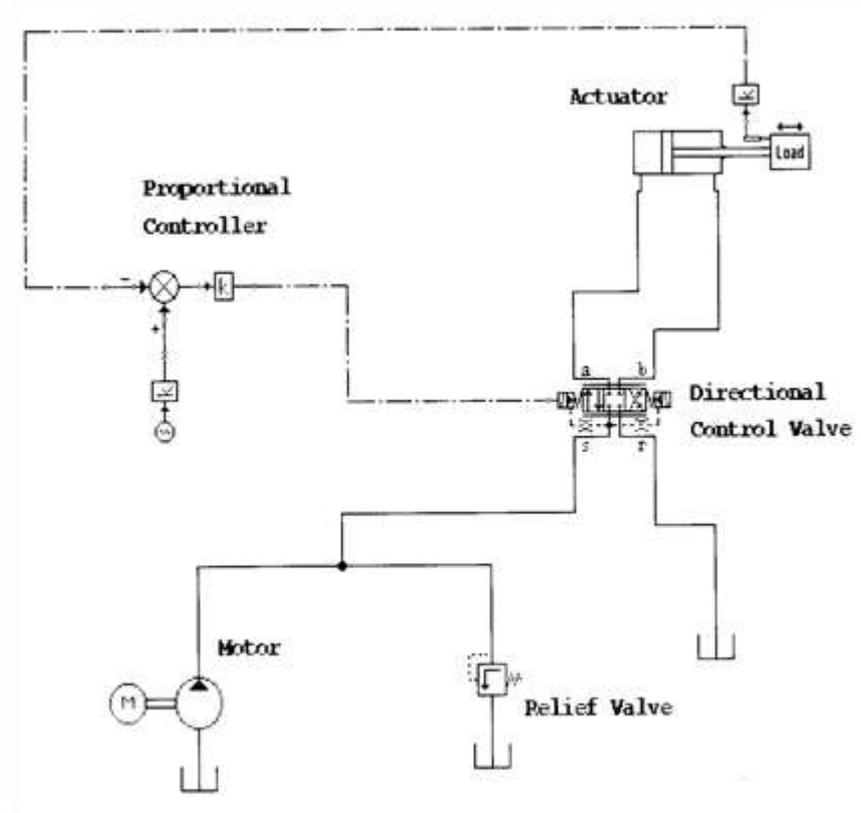

Fig. 4 Linear actuator circuit

## 4.1 Hydrostatic Transmission

In this example the objective of the optimisation is to select capacities for the pump and motor in the hydrostatic transmission so that a given rotational speed is produced by the motor within the first five seconds. A secondary objective is to limit the opening of the relief valve. The load in this circuit consists of a constant load torque of 100Nm. The rotational speed of the prime mover is fixed at 1500 rpm and the cracking pressure of the relief valve is set to 100 bar. The objective function used in this work is given in equation 3.

$$obfn = \left(e\{\omega_m\}\right)^2 \times \left(1 + \frac{Q_{rv}}{Q_p}\right) \tag{3}$$

Squaring the error between the desired and actual operating speeds penalises solutions which have a large error and allows a bandwidth of acceptable error to exist. The percentage of the steady state pump flow that returns to tank through the relief valve is used as a penalty function multiplier in order to force the method to find solutions for which the relief valve does not open.

The pump and motor displacements were constrained between 1 cc/rev and 1000 cc/rev. The initial step size of the search was 64 cc/rev for both variables and the minimum step size was 1 cc/rev. All other parameters are held at constant values. Table 4.1 present results for ten runs of the method on this problem.

Table 1 Parameter sizes for hydrostatic transmission circuit

| RUN | PUMP DISP. (cc/rev) | MOTOR DISP. (cc/rev) | SPEED (r/min) | OBFN | NO. EVALS |
|---|---|---|---|---|---|
| 1 | 125 | 623 | 300.236 | 0.055707 | 622 |
| 2 | 139 | 693 | 300.188 | 0.035425 | 686 |
| 3 | 143 | 713 | 300.176 | 0.031074 | 583 |
| 4 | 128 | 638 | 300.225 | 0.050579 | 666 |
| 5 | 185 | 923 | 300.082 | 0.006667 | 540 |
| 6 | 149 | 713 | 300.160 | 0.025464 | 641 |
| 7 | 160 | 798 | 300.131 | 0.017119 | 582 |
| 8 | 139 | 693 | 300.188 | 0.035425 | 633 |
| 9 | 168 | 838 | 300.115 | 0.013125 | 613 |
| 10 | 100 | 499 | 299.754 | 0.060726 | 606 |

All of these solutions have zero steady state relief valve flow. The graph in Figure 5 shows the speed response of the best circuit, found in the fifth run, for the first five seconds of operation. There are some initial transient components that are damped out so that the motor has a constant nominal operating speed of approximately 300rpm.

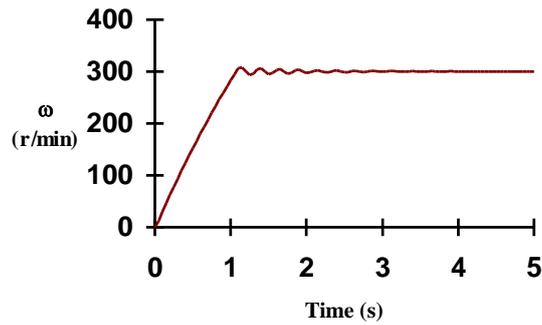

Fig. 5 Speed response of best circuit

This best solution can be compared to a solution found by a hydraulic system designer. In this case the pump and motor displacements are determined by a steady state approach which uses the following equations.

$$T_m = D_m \Delta P_m \eta_{mm} \tag{4}$$

$$Q_m = \frac{D_m \omega_m}{\eta_{vm}} \tag{5}$$

$$Q_p = D_p \omega_p \eta_{vp} \tag{6}$$

Taking efficiencies of 95% for both the pump and motor and assuming the system pressure in the steady state to be 85 bar (i.e. below the relief valve cracking pressure) then the motor displacement required to develop a torque of 100Nm can be calculated to be 78 cc/rev. In steady state the pump and motor flows will be equal and so the required pump displacement is 17.3 cc/rev. In practice, fixed capacity units will not normally be obtainable with the required displacements and as a consequence the motor speed will not be 300rpm. A simulation study undertaken using the calculated displacements and efficiency values given above indicated a problem which would not be apparent to the designer. Although the circuit achieved the required speed within 5 seconds, the system pressure exceeded the relief valve setting during the first 2 seconds of operation. This did not occur with the results obtained by the optimisation method due to the penalty function term included in the objective function. As a result the solutions were forced towards larger unit capacities and lower system pressures. Hence the designer and optimised solutions are not directly comparable as they were obtained using different design objectives.

By changing the form of the objective function to include a minimum pump power term in place of the secondary objective for the relief valve the Tabu search produced solutions close to that obtained using a steady state analysis. For example, one run of the method resulted in a motor displacement of 85 cc/rev and a pump displacement of 17 cc/rev. In this solution the higher motor displacement is due to the use of a lower volumetric efficiency. The terms included in the objective function clearly have a significant effect on the solution obtained and further work is being carried out in this area.

### 4.2 Two Motor Circuit

In this circuit the design objectives are to have two motors operating at two different speeds but driven from the same pump. This is achieved by the use of pressure compensated flow valves (PCFV) and the required circuit is shown in Figure 4.3. A similar problem, though with slightly different parameter settings and objectives, was attempted in previous work (**1**,**2**,**10**) and shown to be a difficult circuit to optimise.

In this example, the design variables selected for optimisation are the displacements of the pump and both motors, along with the nominal flow rates of the two control valves. The pump and motor displacements are constrained between 1cc/rev and 1000 cc/rev. The nominal flow rate settings of the PCFVs are constrained between 10 and 100 L/min. The objective function is given in equation 4.5.

$$obfn = \left(\left|e\{\omega_{m(1)}\}\right| + \left|e\{\omega_{m(2)}\}\right|\right)^2 \times \left(1 + \frac{Q_{rv}}{Q_p}\right) \quad (7)$$

Table 2 shows the parameter sets obtained for ten runs of the problem.

Table 2 Parameter sizes for two motor circuit

| RUN | PUMP (cc/rev) | MOTOR 1 (cc/rev) | MOTOR 2 (cc/rev) | PCF VALVE 1 FLOW RATE (L/min) | PCF VALVE 2 FLOW RATE (L/min) |
|---|---|---|---|---|---|
| 1 | 109 | 883 | 945 | 100.0 | 50.0 |
| 2 | 99 | 916 | 625 | 99.0 | 27.0 |
| 3 | 50 | 510 | 222 | 58.0 | 10.5 |
| 4 | 92 | 723 | 837 | 76.5 | 39.0 |
| 5 | 106 | 712 | 991 | 75.0 | 47.5 |
| 6 | 140 | 526 | 470 | 53.5 | 19.0 |
| 7 | 92 | 643 | 998 | 68.0 | 49.0 |
| 8 | 164 | 927 | 1000 | 100.0 | 48.0 |
| 9 | 108 | 849 | 987 | 95 | 51.5 |
| 10 | 51 | 300 | 656 | 28.0 | 28.5 |

Table 3 shows the objective function value and the number of evaluations for each solution and Table 4 shows the values for the performance characteristics of each of the solutions.

Table 3 Number of evaluations required and objective function values

| RUN | OBFN | NO. EVALS |
|---|---|---|
| 1 | 0.042852 | 1228 |
| 2 | 0.050168 | 1325 |
| 3 | 0.072637 | 1647 |
| 4 | 0.062189 | 1438 |
| 5 | 0.017945 | 1553 |
| 6 | 0.017381 | 1310 |
| 7 | 0.052334 | 1249 |
| 8 | 0.091071 | 1588 |
| 9 | 0.057354 | 1567 |
| 10 | 0.089467 | 1453 |

Table 4 Performance of proposed solutions

| RUN | $\omega_1$ (r/min) | $\omega_2$ (r/min) | RV FLOW (L/min) | PUMP FLOW (L/min) | $\Delta P_1$ (Bar) | $\Delta P_2$ (Bar) |
|---|---|---|---|---|---|---|
| 1 | 120.018 | 59.9761 | 0 | 163.016 | 48.9565 | 54.0694 |
| 2 | 119.979 | 60.0235 | 0 | 147.771 | 88.0096 | 90.9213 |
| 3 | 120.012 | 59.9402 | 0 | 74.8319 | 29.0221 | 27.4343 |
| 4 | 119.968 | 60.0233 | 0 | 137.321 | 86.3007 | 92.6111 |
| 5 | 119.985 | 60.0002 | 12.9587 | 158.202 | 88.1692 | 95.2830 |
| 6 | 119.993 | 60.0001 | 117.3 | 208.944 | 84.9751 | 91.4749 |
| 7 | 119.953 | 60.0024 | 0 | 137.367 | 78.5234 | 86.6793 |
| 8 | 119.993 | 60.0000 | 73.1634 | 244.764 | 90.4574 | 95.4155 |
| 9 | 119.968 | 60.0235 | 0 | 161.457 | 56.4419 | 62.0156 |
| 10 | 120.047 | 60.0401 | 0.39652 | 76.116 | 75.0051 | 93.4306 |

The graph in Figure 6 shows the speed response of the two motors for the solution with the lowest objective function value which was found in the third run. Only the first two seconds are shown as both motors rapidly reach steady state conditions.

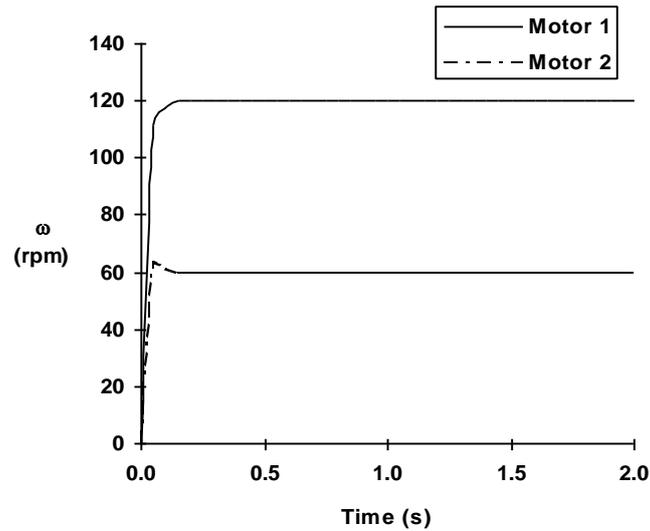

Fig. 6 Speed response

Both motors reach the required steady state speeds. However, considering other factors than pure speed response indicate that this circuit does not exhibit the best characteristics. For example, approximately one quarter of the total pump flow is returning to tank through the relief valve which is an indication that the pump is oversized. However, reducing the pump size does not result in a better solution as the PCFV settings then require altering in order to maintain both an acceptable speed error and a sufficiently high pressure drop across each in order to ensure that they are operating in their linear region. This highlights the complexity of the circuit by showing the interdependencies between each of the system components.

**4.3 Linear Actuator Circuit**

In this final circuit example, the objective is to assign values to circuit components so that the linear actuator produces a desired extend and retract profile when operated under the proportional feedback control. The circuit diagram was shown in Figure 4.

The parameters chosen for inclusion in the optimisation are the pump displacement, the actuator diameter, the actuator stroke, the nominal flow rates of the directional control valve (DCV) and the proportional gain constant. The pump displacement is constrained between 1cc/rev and 500 cc/rev. The actuator diameter is constrained between 30 and 70 mm. The actuator stroke is constrained between 0.01 and 1.0 m. The nominal flow rate settings of the DCV are constrained between 25 and 75 L/min. The proportional gain constant is constrained between 1 and 300.

Because the desired operation of the circuit is a time dependant function, it is important to define an objective function that fully describes this behaviour. The objective function used is given in eqn 4.3 where the subscript $i$ indicates the value of each term at differing time values.

$$obfn = \sum_{i=0}^{n}\left\{\left|x_{di} - x_{ai}\right| \times \left(1 + \frac{Q_{rv}}{Q_p}\right)_i\right\} \quad (8)$$

At each point in the cycle, the position error is calculated and this value penalised by considering the percentage of the pump flow that is flowing through the relief valve.

In order to reduce the total time required for each simulation the number of points where the circuit performance is evaluated has been limited to an interval of 0.1 seconds Table 5 shows the results obtained.

Table 5 Parameter sizes for linear actuator circuit

| RUN | PUMP DISP. (cc/rev) | ACTUATOR DIAMETER (mm) | ACTUATOR STROKE (m) | DCV FLOW RATE (L/min) | GAIN |
|---|---|---|---|---|---|
| 1 | 23 | 52.5 | 0.5 | 69 | 264 |
| 2 | 31 | 60.5 | 0.5 | 75 | 300 |
| 3 | 33 | 63.0 | 0.5 | 75 | 300 |
| 4 | 33 | 63.0 | 0.5 | 75 | 300 |
| 5 | 33 | 63.0 | 0.5 | 75 | 300 |
| 6 | 33 | 63.0 | 0.5 | 75 | 300 |
| 7 | 26 | 56.5 | 0.5 | 73 | 273 |
| 8 | 33 | 63.0 | 0.5 | 75 | 300 |
| 9 | 33 | 63.0 | 0.5 | 75 | 300 |
| 10 | 33 | 63.0 | 0.5 | 75 | 300 |

It can be seen that the parameter values obtained in each run are much more similar than in the previous examples. This implies that the definition of the objective function is much more rigorous and only a small number of locally optimum solutions exist. Table 6 shows the objective function value for each solution along with the number of evaluations that was required.

Table 6 Number of evaluations required and objective function values

| RUN | NO. EVALS | OBFN |
|---|---|---|
| 1 | 1307 | 0.032576 |
| 2 | 1329 | 0.055856 |
| 3 | 1378 | 0.061436 |
| 4 | 1117 | 0.061436 |
| 5 | 1240 | 0.061436 |
| 6 | 1390 | 0.061436 |
| 7 | 1608 | 0.053672 |
| 8 | 1745 | 0.061436 |
| 9 | 1121 | 0.061436 |
| 10 | 1222 | 0.061436 |

Figure 7 shows the output motion of the actuator for the best solution, found in the first run, in comparison to the desired profile. In this simulation the number of points evaluated has been increased by changing the time step interval to 0.05 seconds. This compares the simulation used in the optimisation process where the time step interval was 0.1 seconds.

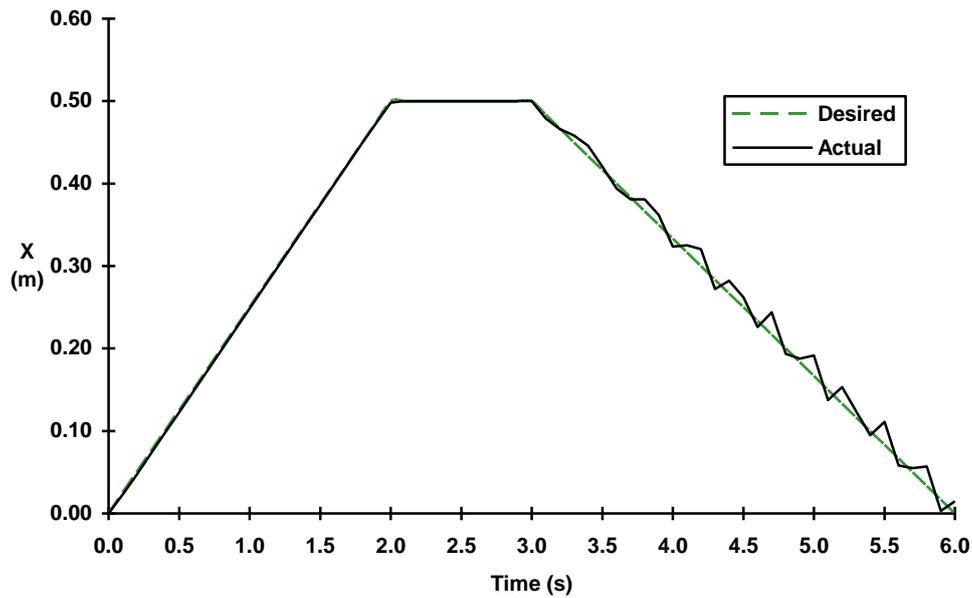

Fig. 7 Actual and desired motion profile of best solution

Despite the relatively low objective function value it can be seen that considerable errors exist due to the stability of the system. However, these errors are unlikely to be as apparent during the optimisation process due to the reduced number of evaluation points.

However, considering the solution which was found most consistently it can be seen that the method is able to locate viable solutions as this solution locates a relatively low error despite the reduced number of evaluation points. The response of this circuit is shown in Figure 8.

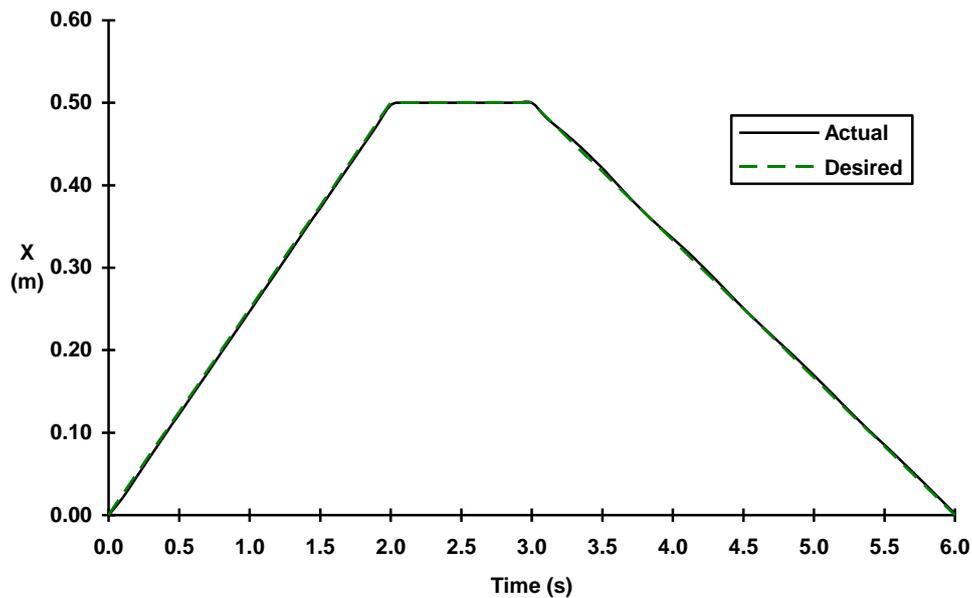

Fig. 8 Actual and desired motion profile of common solution

## 5  DISCUSSION

The Tabu search method described in this paper offers an aggressive technique that is applicable to many optimisation problems, including hydraulic circuit design. The results presented for the three different circuits illustrate the effectiveness of the method. However, they also show that considerable work is still required to enable real design objectives to be included into a purely numerical objective function.

Considering the approach used by a human design engineer provides considerable insight and potential to improve the method. In the example used in this paper, that of steady state calculations, it can be seen that the approach does not lead to ideal solutions. However, considering the use of simulation by a human designer on the more complex circuits provides more useful information. In all cases the engineer used simple calculations to obtain a first estimate and then refined this solution to obtain the desired performance. This estimation of an initial start point combined with a reduction in the size of the search space is likely to lead to a much more efficient optimisation.

Several other points were apparent. The first of these was the tendency of the human designer to increase the cracking pressure of the relief valve in order to obtain acceptable performance. It may be beneficial in future work to include the relief valve cracking pressure as a design variable in order to achieve greater flexibility of solution. For circuits where low system pressure is essential a penalty function could be used to limit the cracking pressure of the relief valve.

A second point was the tendency of the designer to use as small a pump as possible, so limiting the power consumption of the circuit. An additional penalty function based on pump size could be used to force the method towards solutions with low power consumption.

However, the inclusion of large numbers of penalty functions is likely to lead to a reduction in performance of the search method due to the interdependence of the system components. This will be realised as large numbers of local optima in the objective function. One possible method for eliminating such large numbers of penalty functions would be through the implementation of a multi-objective search.

## 6  CONCLUSIONS

This paper has shown that Tabu search techniques can be successfully applied to hydraulic circuit optimisation and are considerably more efficient than many other optimisation methods. This increase in efficiency is due to the reduced number of objective function evaluations required to find a solution. This is particularly true for when design parameter values are restricted to sensible, discrete values.

However, whilst the Tabu search locates acceptable solutions to a given problem it appears to lack the consistency that is required if a single run is expected to locate the optimal solution. To a certain extent this is due to poor definition of objective function. This poor definition is mainly due to the necessity of incorporating more than a single design objective into a single function. In this paper the main component is penalised by multiplying it by a penalty ratio. This approach eliminates the need for the selection of weighting factors used in the more common weighted sum approach but still lacks the power of a true multi-objective search.

Further work is required to deal with multiple objectives in fluid power system design. Initially, such work should consider the use of domination to determine the pareto-optimal solutions found by the method. However, even in this approach it will be necessary to consider the definition of allowable moves from non-dominated to dominated solutions.

Once this issue has been dealt with enhancements to the multi-objective approach could be considered which may incorporate aspects of fuzzy logic to determine the degree of domination of similar solutions.

# 7  ACKNOWLEDGMENTS

The research in this paper is funded by the Engineering and Physical Sciences Research Council under grant no. 86796 and this support is gratefully acknowledged.